\documentclass[10pt,a4paper]{article}

\usepackage[utf8]{inputenc}
\usepackage[english]{babel}
\usepackage{amsmath}
\usepackage{amsfonts}
\usepackage{amssymb}
\usepackage{makeidx}
\usepackage{graphicx}
\newcommand{\intt}{\int_0 ^t}
\newcommand{\intT}{\int_0 ^T}

 \newcommand{\intzT}{\int_\zt^T}
  \newcommand{\intzt}{\int_\zt^t}
\newcommand{\Xz}{X_{\zt }}
\newcommand{\zt}{\tau}
\newcommand{\ZR}{\rangle} 
\newcommand{\ZL}{\langle}

\newcommand{\zdia}{~~\rule{1mm}{2mm}\par\medskip}  

\newcommand{\ZIN}{\infty}

\newcommand{\zzr}{{\rm I\hskip-2.1pt R}}

\newcommand{\ZBI}{\bibitem}
\newcommand{\ZCD}{(\cdot)}
\newcommand{\ZD}{\;\mbox{\rm d}}

\newcommand{\ZLA}{\label}

\newtheorem{Theorem}{Theorem}

\newtheorem{Remark}[Theorem]{Remark}
\newtheorem{Definition}[Theorem]{Definition} 
\newtheorem{Example}[Theorem]{Example}

\author{
L. Pandolfi\thanks{Dipartimento di Scienze Matematiche ``Giuseppe Luigi Lagrange'', Politecnico di Torino, Corso Duca degli Abruzzi 24, 10129 Torino, Italy (luciano.pandolfi@polito.it)}
}
\title{The quadratic regulator problem and the Riccati equation for a process governed by a linear Volterra integrodifferential equations\thanks{
This papers fits into the research program of the GNAMPA-INDAM and has been written in the framework of the   ``Groupement de Recherche en Contr\^ole des EDP entre la France et l'Italie (CONEDP-CNRS)''.}}
\begin{document}
\maketitle

{\bf \underline{Abstract}:} In this paper we study the quadratic regulator problem for a process governed by a Volterra integral equation in $\zzr^n$. Our main goal is the proof that it is possible to associate a Riccati differential equation to this quadratic control problem, which leads to the feedback form of the optimal control. This is in contrast with previous papers on the subject, which confine themselves to study the Fredholm integral equation which is solved by the optimal control.



\section{Introduction}

The quadratic regulator problem for control processes regulated by linear differential equations both in finite and infinite dimensional spaces has been at the center of control theory   at least during the last eighty years, after the proof that  the synthesis of dissipative systems amounts to the study of a (singular) quadratic control problem (see~\cite{Brune}). In this period, the theory reached a high level of maturity and the 
monographs~\cite{Bittanti,lasieckaTriggENcicl} contain  the crucial ideas used in the study of the   quadratic regulator problems for lumped and distributed systems (see~\cite{BucciPANDO1,BucciPANDO2,PandSing1,Pandsing2,Pandsing3} for the singular quadratic 
regulator problem for distributed systems). 

In recent times, the study of controllability of systems described by Voterra integrodifferential equations (in Hilbert spaces) has been stimulated by several applications (see~\cite{Pandlibro}) while the theory of the quadratic regulator problem for these   systems is still at a basic level. In essence, we can cite only the paper~\cite{Pritch} and some applications of the results in this paper, see for example~\cite{HUANGliWANG}. In these papers, the authors study a standard regulator problem for a system governed by a Volterra integral equation (in a Hilbert space and with bounded operators. The paper~\cite{HUANGliWANG} and some other applications of the results in~~\cite{Pritch} studies a stochastic system) and the synthesis of the optimal control is given by relying on the usual variational approach and Fredholm integral equation for the optimal control.
The authors of these papers do not develop a Riccati differential equation   and this is our goal here. In order to avoid the technicalities   inevitably introduce by the presence of unbounded operators which are introduced by the action of boundary controls,  we confine ourselves to study Volterra integral equations in $\zzr^n$.
 
The control problem we consider is described by
\begin{equation}
\ZLA{eq:Volte}
x'=\intt N(t-s)x(s)\ZD s+ Bu(t)\,,\qquad x(0)=x_0 
\end{equation} 
where
$x\in \zzr^n$, $u\in \zzr^m$, $B$ is a constant $n\times m$ matrix and $N(t)$ is a continuous $n\times n$ matrix (extension to $B=B(t)$ and $N=N(t,s)$ is simple). 
Our goal is the study of the minimization of the standard quadratic cost
\begin{equation}
\ZLA{eq:costoAZERO}
\intT\left [ x^*(t) Qx(t) +|u(t)|^2\right ]\ZD t+ x^*(T)Q_0x(T) 
\end{equation}
where $Q=Q^*\geq 0$, $Q_0=Q_0^*\geq 0$.  
 
Existence of a unique optimal control in $L^2(0,T;\zzr^m)$ for every fixed $x_0\in \zzr^n$ is obvious.

The plan of the paper is as follows: in order to derive a Riccati differential equation, we need a suitable ``state space'' in which our system evolves. In fact, a Volterra integral equation is a semigroup system in a suitable infinite dimensional space (see~\cite[Ch.~6]{Nagel}) and we could relay on this representation of the Volterra equation to derive a theory of the Riccati equation in a standard way but the shortcoming   is that the ``state space'' is 
$  \zzr^n\times L^2(0,+\ZIN ;\zzr^n)$ and the Riccati differential equation so obtained should be solved in a space with infinite memory, even if the process is considered on a finite time interval $[0,T]$. We wish a ``Riccati differential equation'' in a space which has a ``short memory'', say of duration at most $T$, as required by the optimization problem. So, we need the introduction of a different ``state space approach'' to Eq.~(\ref {eq:Volte}). This is done in Sect.~\ref{sec:SateDyn} where, using dynamic programming, we prove that the minimum of the cost is a quadratic form which satisfy a 
(suitable version) of the Linear Operator Inequality {\bf (LOI)}.

Differentiability properties of the cost are studied  in section~\ref{sect:DIFFEpropri} (using a variational approach to the optimal control related to the arguments in~\cite{Pritch}). The regularity properties we obtain finally allows us to write explicitly a system of partial differential equations (with a quadratic nonlinearity) on $[0,T]$, which is the version of the Riccati differential equations for our system.

We believe that the introduction of the state space in Sect.~\ref{sec:SateDyn} is a novelty of this paper.

\section{\ZLA{sec:SateDyn}The state of the Volterra integral equation, and the {\bf (LOI)}}

  According to the general definition in~\cite{Kalman}), the state at time $\zt $ is the   information at time $\zt $ needed to uniquely solve the equation for $t>\zt $ (assuming the control is known for $t>\zt $).

It is clear that if $\zt =0$ then the sole vector $x_0$ is sufficient to solve   equation~(\ref{eq:Volte}) in the future, and the state space at $\zt =0$ is $\zzr^n$. Things are different if we solve the equation till time $\zt $ and we want to solve it in the future. In this case, Eq.~(\ref{eq:Volte}) for $t>\zt $ takes the form
\begin{equation}\ZLA{eq:volteat0}
x'=\int_{\zt }^t N(t-s)x(s)\ZD s+B u(t)+ \int_0^{\zt } N(t-s)x(s)\ZD s\,.
\end{equation}
In order to solve this equation for $t>\zt $ we must know the pair\footnote{Remark on the notation: $x_{\zt }=x_{\zt }(s)$ is a function on $(0,\zt )$ while $X_{\zt }$ (upper case letter) is the pair $(x(\zt ),x_{\zt })$.}
 $X_{\zt }=\left (x(\zt ),x_{\zt }(\cdot) \right )$ where $x_{\zt }(s) =x(s)$, $s\in (0,\zt )$.

Note that in order to uniquely solve~(\ref{eq:volteat0}), $x_{\zt }(\cdot)$ needs not be a segment of previously computed trajectory. It can be an ``arbitrary'' function.
This observation suggests the definition of the following state space at time $\zt $: 
\[
M^2_{\zt }=\zzr^n\times L^2(0,\zt ;\zzr^n)
\]
(to be compare with the state space of  differential equations with a fixed delay $h$ which is $\zzr^n\times L^2(-h,0;\zzr^n)$).

  Eq.~(\ref{eq:volteat0}) defines, for every fixed $u$ and $\zt_1>\zt $, a solution map from $M^2_{\zt }$ to $M^2_{\zt_1}$ which is affine linear and continuous.
An explicit expression of this map can be obtained easily. Let us fix an initial time $\zt\geq 0$. Let  $t\geq \zt$ and let $Z(t,\zt )$ be the $n\times n$ matrix solution of
\begin{equation}\ZLA{eq:diZgrande}
\frac{\ZD}{\ZD t}Z(t,\zt )=\int_{\zt }^t Z(\xi,\zt)N(t-\xi) \ZD \xi\,,\quad Z(\zt,\zt)=I \,.
\end{equation}
 Then,
 
\begin{equation}\ZLA{eq:evoluzione}
x(t)=Z(t,\zt ) \hat x  +\int_0 ^{\zt } Y(t,s;\zt ) \tilde x(s)\ZD s+\int _{\zt }^t Z(t-r+\zt,\zt)B u(r)\ZD r 
\end{equation}
where
\[
Y(t,s;\zt )=\int_{\zt }^t Z(t-\xi+\zt,\zt)N(\xi-s)\ZD \xi\,.
\]

This way,  for every $\zt_{1}>\zt $ we define two linear continuous transformations: $E(\zt_1;\zt )$ from $M^2_{\zt }$ to $M^2_{\zt_1}$ (when $u=0$) and $\Lambda(\zt_1;\zt )$ from $L^2(\zt ,\zt_1;\zzr^m)$ to $M^2_{\zt_1}$ (when $X_{\zt }=0$), as follows:
\[
E(\zt_1;\zt )(\hat x,\tilde x(\cdot))=(x(\zt_1),y)\qquad y=\left\{
\begin{array}{lll}
x(t)\ \mbox{given by~(\ref{eq:evoluzione})}&{\rm if}& \zt <t<t_1\\
\tilde x(t) &{\rm if} & t\in (0,\zt)\,.
\end{array}
\right.
\]
The operator $\Lambda(\zt_1;\zt )$ is defined by the same formula as $E(\zt_1;\zt )$, but when $X_{\zt }=0$ and $u\neq 0$.

The evolution of the system is describe by the operator
\begin{equation}\ZLA{eq:evoluSISTE}
 E(t_1;\zt )X_{\zt }+\Lambda(t_1;\zt )u\,.
\end{equation}

The evolutionary properties of this operator  follow  from the unicity of solutions of the Volterra integral equation. Let us consider Eq.~(\ref{eq:volteat0}) on $[\zt,T]$ with initial condition 
$(\hat x,\tilde x(\cdot))$, whose solution is given by~(\ref{eq:evoluzione}). Let $\zt_1\in (\zt,T)$ and let us consider Eq.~(\ref{eq:volteat0}) on $[\zt_1,T]$ but with initial condition $\left (x(\zt_1),x_{\zt_1}\right )$.
Eq.~(\ref{eq:volteat0}) on $[\zt_1,T]$ and this initial condition takes the form
\[
x'(t)=\int _{\zt_1} ^t N(t-s) x(s)\ZD s+Bu(t)+\int_0^{\zt_1}N(t-s) x_{\zt_1}(s)\ZD s\,,\qquad x(\zt_1^+)=x(\zt_1^-)
\]
and so, on $[\zt_1,T]$ we have
\[
x'(t)=Z(t,\zt_1)x(\zt_1^-) +\int_0^{\zt_1} Y(t,s;\zt_1) x_{\zt_1}(s)\ZD s+\int _{\zt_1}^t Z(t-s-\zt_1,\zt_1)Bu(s)\ZD s\,.
\]
Unicity of the solutions of the Volterra integral equation shows that, for $t\in (\zt_1,T]$ the following equality holds 
\[
E(t,\zt)\left (\hat x,\tilde x\right )+\Lambda(t;\zt )u=E(t,\zt_1)\left [ E(\zt_1,\zt)\left (\hat x,\tilde x\right )+\Lambda(\zt_1;\zt )u\right ]+ \Lambda(t;t_1)u\,.
\]

\begin{Remark}\ZLA{rema:DerivZistiniz}
{\rm
The solution $Z(t,\zt)$ of Eq.~(\ref{eq:diZgrande}) solves the following Volterra integral equation on $[\zt,T]$:
\[
Z(t)=1+\int_\zt^t Z(\xi)M(t-\xi)\ZD\xi\,,\qquad M(t)=\intt N(s)\ZD s\,.
\]
The usual Picard iteration gives
\begin{align*}
Z(t,\zt)&=1+\int_\zt^t M(t-\xi)\ZD\xi+\int_\zt^t \int_\zt^\xi M(\xi-\xi_1)\ZD\xi_1 M(t-\xi)\ZD\xi+\cdots=\\
&= 1+\int_\zt^t M(t-\xi)\ZD\xi+\int_\zt^t \int_0^{t-s} M(t-s-r)M(r)\ZD r\ZD s+\cdots 
\end{align*}
 The properties of these integrals is that, once exchanged, we have
 \[
Z(t,\zt)=1+\int_\zt^t H(t-s)\ZD s 
 \]
 where $H(t)$ does not depend on $\zt$ and it is differentiable. It follows that \emph{the function $(\zt,t)\mapsto Z(t,\zt)$ is continuously differentiable on $0<\zt<t<T$ and the derivative has continuous extension to $0\leq \zt\leq t\leq T$.}\zdia
}
\end{Remark}

Now we begin our study of the quadratic regulator problem and of the Riccati equation.

One of the possible ways to derive an expression of the optimal control and possibly  a Riccati differential equation for the quadratic regulator problem is via dynamic programming. We follow this way.  For every fixed $\zt <T$ we introduce
\[
J_{\zt }\left (X_{\zt },u\right )=\intzT \left [ x^*(t)Qx(t)  +|u(t)|^2\right ]\ZD t  +x^*(T)Q_0x(T)
\]
where $x(t)$ is the solution of~(\ref{eq:volteat0}) (given by~(\ref{eq:evoluzione})) and we define
\begin{equation}\ZLA{eq:identiLOI}
W(\zt ;\Xz)=\min _{u\in L^2(\zt ,T;\zzr^m)} J_{\zt }\left (\Xz,u\right )\,.
\end{equation}
Existence of the minimum is obvious and we denote $u^+(t)=u^+(t;\zt ,\Xz)$ the optimal control.
 The corresponding solution is denoted $x^+(t)=x^+(t;\zt ,\Xz)$ while we put $X^+_{t }=\left (x^+(t ),x^+_{t }\ZCD\right )$. 

Let us fix any ${\zt_1} \in (\zt ,T)$ and let $u(t)=u^1(t)$ if $t\in (\zt ,{\zt_1} )$, $u(t)=u^2(t)$ if $t\in ({\zt_1} ,T)$, while
\[
X^1_t=E(t,\zt )\Xz+\Lambda(t,\zt )u^1\quad  t\in[\zt ,{\zt_1} ]\,,\quad 
X^2_t=E(t,{\zt_1} )X^1_{{\zt_1}} +\Lambda(t,{\zt_1} )u^2\quad t\in[{\zt_1} ,T]\,.
\]
We noted that $X (t;\zt ,X_{\zt })$ given by~(\ref{eq:evoluSISTE}) on $[\zt ,T]$ is equal to $X^1_t$ on $[\zt ,{\zt_1} ]$ and to $X^2_t$ on $[{\zt_1} ,T]$.

Let $x^i$ be the $\zzr^n$ component of $X^i$. Then, for every $u$ we have (we use the crochet to denote the inner product instead of the more cumberstome notation $\left  (x^1(t)\right )^*Q x_1(t)$)
\begin{equation}\ZLA{eq:PREloi}
W(\zt ,\Xz)\leq \int_{\zt }^{\zt_1}  \left [\ZL Qx^1(t),x^1(t)\ZR +|u^1(t)|^2\right ]\ZD t+J_{\zt_1} \left ( X^1_{\zt_1} ,u^2\right )\,.
\end{equation}
This inequality holds for every $u^1$ and $u^2$ and equality holds when $u^1$ and $u^2$ are restrictions of the optimal control $u^+$.

 We keep $u^1$ fixed and 
we compute the minumum of the right hand side respect to $u^2$. We get the Linear Operator Inequality {\bf (LOI)}:
\begin{equation}\ZLA{eqDiseqLOI}
W(\zt ,X_{\zt })\leq \int _{\zt }^{\zt_1}  \left [\ZL Qx^1(t),x^1(t)\ZR+|u^1(t)|^2\right ]\ZD t+W\left ({\zt_1} ,X^1_{\zt_1} \right )\,. 
\end{equation}
This inequality holds for every control $u\in L^2(\zt ,{\zt_1} ;\zzr^n)$.
Let in particular $u^1$ be the restriction to $(\zt ,{\zt_1} )$ of $u^+(\cdot)=u^+(\cdot;\zt ,X_{\zt })$. Inequality~(\ref{eq:PREloi}) shows that the minimum of 
$J_{\zt_1} \left ( X^1_{\zt_1} ,u^2\right )$ cannot be strictly less then $J_{\zt_1} \left ( X^1_{\zt_1} ,u^+\right )$, i.e. 
 the optimal control of the cost $J_{\zt_1} \left ( X^1_{\zt_1} ,u^2\right )$ is the restriction to  $({\zt_1} ,T)$ of $u^+(t)$, the optimal control of $J_{\zt }\left (X_{\zt },u\right )$. 

Equality holds in~(\ref{eqDiseqLOI}) if $u^1=u^+$.

 In conclusion, we divide with ${\zt_1-\zt} $ (which is positive) and we find the following inequality, \emph{which holds with equality if $u=u^+$:}
\[
\frac{1}{{\zt_1-\zt} }\left [W\left ({\zt_1} ;X^1_{\zt_1} \right )-W\left ( \zt ;X_{\zt  }\right )\right ]\geq -\frac{1}{{\zt_1-\zt} } \int _{\zt }^{\zt_1}  \left [
\ZL Qx^1(t),x^1(t)\ZR+|u(t)|^2
\right ]\ZD t\,.
\]
  So,   the following inequality   holds when $\zt $ is a Lebesgue point of $u(t)$ (every $\zt $ if  $u$ is continuous):
\begin{equation}\ZLA{eq:PrimaFormINEqEQ}
\lim\inf _{{\zt_1} \to \zt ^+}\frac{1}{{\zt_1-\zt} }\left [W\left ({\zt_1} ;X^1_{\zt_1} \right )-W\left ( \zt ;X_{\zt  }\right )\right ]\geq -\left [ \ZL Q x(\zt ),x(\zt )\ZR+|u(\zt )|^2\right ]\,.
\end{equation}
Equality holds if $u=u^+$ and $\zt $ is a Lebesgue point of $u^+$ and in this case we can even replace $\liminf$ with $\lim$, i.e. $W\left ({\zt_1} ;X^+_{\zt_1} \right )$ is differentiable if $\zt $ is a Lebesgue points   of $u^+$.

The previous argument can be repeated for every $\zt $ so that the previous inequalities/equalities holds $a.e.$ on $[0,T]$ and   we might even replace $\zt $ with the generic notation $t$.

\begin{Remark}
{\rm 
If it happens that $\ker N(t)= S$, a subspace of $\zzr^n$, we might also consider as   the second component of the ``state'' $X_{\zt }$ the projection of $\tilde x$ on (any fixed) complement of $S$, similar to the theory developed in~\cite{DelfourMANITIUS,FABRIZIOPATA}. We dont't pursue this approach here.\zdia
}
\end{Remark}



\section{\ZLA{sect:DIFFEpropri}The regularity properties of the value function, the synthesis of the optimal control  and the Riccati equation}

We  prove that $W$ is a continuous quadratic form with smooth coefficients and we prove that $u^+(t)$ is continuous (so that every time $t$ is a Lebesgue point of $u^+(t)$). We arrive at this result via the variational characterization of the optimal pair $(u^+,x^+)$ ($x^+$ is the $\zzr^n$-component of $X^+$) in the style 
of~\cite{Pritch}. The standard perturbation approach gives a representation of the optimal control (and a definition of the adjoint state $p(t)$): 
\begin{equation}\ZLA{eq:defiOTTIMcontroVariaz}
u^+(t)=-B^*\left [\int_t^T Z^*(s-t+\zt,\zt)Q x^+ (r)\ZD r  +Z^*(T-t+\zt,\zt)Q_0x^+(T)          \right ]=-B^* p(t)
\end{equation}
where $p$, the function in the bracket,  solves the adjoint   equation
\begin{equation}
\ZLA{eq:aggiunta}
p'(t)=-Qx^+(t)-\int_t^T N^*(s-t)p(s)\ZD s\,,\qquad p(T)=Q_0 x^+(T)\,.
\end{equation}

Note that  $p$ depends on $\zt$ and that Eq.~(\ref{eq:aggiunta}) has to be solved (backward) on the interval $[\zt,T]$.

The simplest way to derive the differential equation~(\ref{eq:aggiunta}) is  to note that  the function $q(t)=p(T-t)$ is given by
\begin{align*}
q(t)&=\int_{T-t}^T Z^*(s-T+\zt+t,t)Qx^+(s)\ZD s+Z^*(t+\zt,\zt)Q_0x^+(T)=\\
&=\intt Z^*(t-r+\zt,\zt)Qx^+(T-r)\ZD r+Z^*(t+\zt,\zt)Q_0x^+(T)\,.
\end{align*}
Comparison with~(\ref{eq:evoluzione}) shows that $q(t)$
solves
\[
q'(t)=\intt N^*(t-s)q(s)\ZD s+Qx^+(T-t)\,,\qquad q(0)=Q_0x^+(T)
\]
from which the equation of $p(t)$ is easily obtained.

We recapitulate: the   equations which characterize $(x^+,u^+)$ when the initial time is $\zt $ and $X_{\zt }=\left (\hat x,\tilde x(\cdot)\right )$ is the following system  of equations on the interval $[\zt,T]$:
\begin{equation}
\begin{array}{ll}
\ZLA{eq:coppiaottima}
x'=\intzt N(t-s)x(s)\ZD s-BB^*p(t)+\int_0^{\zt }N(t-s) \tilde x(s)\ZD s\,, & x(\zt )=\hat x\\[2mm]
p'(t)=-Qx(t)-\int_t^T N^*(s-t) p(s)\ZD s\,, & p(T)=Q_0 x(T)\\[2mm]
u^+(t)=-B^*p(t)\,.
\end{array}
\end{equation}

We replace $u^+(t)=u^+(t;\zt ,X_{\zt }) $ in~(\ref{eq:evoluzione}). The solution is $x^+(t)$. Then we replace the resulting expression in~(\ref{eq:defiOTTIMcontroVariaz}). We get the  Fredholm integral equation for $u^+(t)$:  
\begin{align*}
 u^+(t)&+B^* Z^*(T-t+\zt,\zt) Q_0\int _{\zt }^T Z(T-r+\zt,\zt)B u^+(r)\ZD r+\\
 &+B^*\int_t^T Z^*(s-t+\zt,\zt)Q
 \int_\zt^s Z(s-r+\zt,\zt)Bu^+(r)\ZD r\ZD s=\\
&= -B^*\left [
Z^* (T-t+\zt,\zt)Q_0F(T,\zt)+ 
 \int_t^T Z^* (s-t+\zt,\zt)QF(s,\zt)\ZD s\right ]
\end{align*}
where
\[
F(t,\zt)=Z(t,\zt)\hat x+\int_0^\zt Y(t,s;\zt) \tilde x(s)\ZD s\,.
\]
This Fredholm integral equation has to be solved on $[\zt,T]$.

By solving the Fredholm integral equation we find an expression for  $u^+(t)$, of the following form:
\begin{equation}\ZLA{eq:opeloopcontrol}
u^+(t)=u^+(t;\zt ,X_{\zt }) = \Phi_1(t,\zt ) \hat x+\int_0 ^{\zt } \Phi_2(t,s;\zt ) \tilde x(s)\ZD s\,,\quad t\geq \zt 
\end{equation}
and so also

\begin{equation} \ZLA{eq:openloopSTATE}
x^+(t)=x^+(t;\zt ,X_{\zt }) =Z_1(t,\zt )\hat x+\int_0 ^{\zt } Z_2(t,r;\zt ) \tilde x(r)\ZD r\,, \quad t\geq \zt \,.
\end{equation}
The explicit form of the matrices $\Phi_1(t,\zt )$, $\Phi_2(t,s;\zt )$, $Z_1(t,\zt )$, $Z_2(t,r;\zt )$    (easily derived using the resolvent operator of the Fredholm integral equation) is not needed. The important fact is that \emph{these matrices have continuous partial derivative respect to their arguments $t$, $s$ and $\zt $.  In particular, $u^+(t)=u^+(t;\zt ,X_{\zt })$ is a continuous function of $t$ for $t\geq \zt $. The derivative has continuous extensions to $s=\zt$ and to $t=\zt$.} Differentiability respect to $\zt$ follows from Remark~\ref{rema:DerivZistiniz}.

We replace~(\ref{eq:opeloopcontrol}) and~(\ref{eq:openloopSTATE}) in~(\ref{eq:identiLOI}) and we get 
\begin{align}
\nonumber W(\zt ;X_{\zt })&=\int_{\zt } ^T \left |
Q^{1/2}Z_1(s,\zt )\hat x +Q^{1/2}\int _0^{\zt } Z_2(s,r;\zt ) \tilde x(r)\ZD r
\right |^2\ZD s+\\
\ZLA{eq:FinaPerInfoootimo}&+ \int_{\zt }^T \left |
\Phi_1(s;\zt )x_0+\int_0 ^{\zt } \Phi_2(s,r;\zt )\tilde x(r)\ZD r
\right |^2\ZD s\,.
\end{align}
This equality shows that  $X_{\zt} \mapsto  W(\zt ,X_{\zt })$ is a continuous quadratic form of $X_{\zt }\in M_{\zt }$.

We use dynamic programming again, in particular the fact that $u^+(\cdot;{\zt_1} ,X^+_{{\zt_1} })$ is the restriction to $[{\zt_1} ,T]$ of $u^+(\cdot;\zt ,X_{\zt })$. Hence, for every ${\zt_1} \geq \zt $ we have
 \begin{align}
\nonumber W({\zt_1}  ;X^+_{{\zt_1} })&=\int_{{\zt_1}  } ^T \left |
Q^{1/2}Z_1(s,{\zt_1}  )  x^+({\zt_1} ) +Q^{1/2}\int _0^{{\zt_1} } Z_2(s,r;{\zt_1}  )   x^+(r)\ZD r
\right |^2\ZD s+\\
\ZLA{eq:FinaPerInfoootimo}&+ \int_{{\zt_1}   }^T \left |
\Phi_1(s;{\zt_1} )x^+({\zt_1} )+\int_0 ^{{\zt_1}  } \Phi_2(s,r;t{\zt_1}  )  x^+(r)\ZD r
\right |^2\ZD s\,.
\end{align}

 We simplify the notations: from now on we drop the ${}^+$ and we replace ${\zt_1} $ with $t$ but we must recall that we are computing for $t\geq \zt $ and, when we use equality in~(\ref{eqDiseqLOI}),  on the optimal evolution.
 
 By expanding the squares  we see that $W({\zt_1}  ;X_{{\zt_1}  })$ has the following general form:
\begin{align}\nonumber
W(t  ;X_{t  })&= x^* (t )  P_0( t)x ( t)+ x^*( t )  \intt  P_1( t,s)  x (s)\ZD s+\\
&\ZLA{eq:ExprreDELLAformAquaDRATI}
+\left [\intt P_1(t,s)  x (s)\ZD s\right ]^* x (t)+\intt\intt   x^* (r)  K(t,\xi,r)   x (\xi)\ZD \xi\ZD r\,.
\end{align}
For example,   
\[
P_0(t)= \int_t^T \left [Z_1^*(s,t)QZ_1(s,t)+\Phi_1^*(s,t)\Phi_1(s,t)\right ]\ZD s\,.
\]
Note that $P_0(t)$ is a selfadjoint differentiable matrix.

Now we consider the matrix $ K(t,\xi,r)$.
We consider the contribution of the first line in~(\ref{eq:ExprreDELLAformAquaDRATI}) (the contribution of the second line is similar). Exchanging the order of integration and the names of the variables of integration, we see that
\begin{align*}
& \intt  x^*(r) K(t,\xi,r)   x(\xi)\ZD \xi\ZD r= \intt\intt x^*(r)\left [
\int_t^T Z_2^*(s,r,t)QZ_2(s,\xi,t)\ZD s
\right ]x(\xi)\ZD r\ZD\xi=\\
&=\intt\intt x^*(\xi)\left [
\int_t^T Z_2^*(s,\xi,t)QZ_2(s,r,t)\ZD s
\right ]x(r)\ZD\xi\ZD r=\int_0^t\int_0^t 
x^*(\xi) K^*(t,r,\xi) x(r)\ZD\xi\ZD r
\end{align*}
so that we have
\[
K(t,\xi,r)=K^*(t,r,\xi) 
\]
and this matrix function is differentiable respect to its arguments $t$, $r$ and $\xi$.

Analogously we see differentiability of $P_1(t,s)$.

We whish a differential equations for the matrix functions $P_0(t)$, $P_1(t,s)$, $K(t,s,r)$. In order to achieve this goal, we
  compute the right derivative of $W(t;X _t)$ (and any continuous control) for $t>\zt $ and we use inequality~(\ref{eq:PrimaFormINEqEQ}). We use explicitly that equality holds in~(\ref{eq:PrimaFormINEqEQ}) when the derivative is computed along an optimal evolution. 
  
\subsection{The Riccati equation}

In order to derive a set of differential equations for the matrices $P_0(t)$, $P_1(t,s)$, $K(t,\xi,r)$ we proceed as follows: we fix (any) $\zt\in [0,T]$ and the initial condition $X_{\zt}=(\hat x,\tilde x(\cdot))$. We consider~(\ref{eq:ExprreDELLAformAquaDRATI}) with any continuous control $u(t)$ on $[\zt,T]$ (the corresponding solution of the Volterra equation is $x(t)$). We consider the quadratic form $W$ with the control $u(t)$ and the corresponding solution $X_t$ given in in~(\ref{eq:ExprreDELLAformAquaDRATI}). In this form we separate the contribution of the functions on $(0,\zt)$ and the contribution on $[\zt,t]$. For example $x^*(t)P_0(t)x(t)$ remains unchanged while $x^*(t)\intt P_1(t,\xi)x(\xi)\ZD \xi$ is written as 
\[
x^*(t)\intt P_1(t,\xi)x(\xi)\ZD \xi=x^*(t)\int_0^\zt  P_1(t,s)\tilde x(s)\ZD s+x^*(t)\int_\zt^t P_1(t,s) x(s)\ZD s\,.
\]
The other addenda are treated analogously.

We obtain a function of $t$ which is continuously differentiable. Its derivative at $t=\zt$ is the left hand side of~(\ref{eq:PrimaFormINEqEQ}) and so it satisfies the inequality~(\ref{eq:PrimaFormINEqEQ}), with equality if it happens that we compute with $u=u^+$. So,  the function of $u\in\zzr^m$
\[
u\mapsto \left [\frac{\ZD}{\ZD t}W(\zt  ;X _{\zt  })+u^*(\zt)u(\zt)\right ]
\]
reaches a minimum at $u=u^+_\zt$. Note that $\zt\in [0,T]$ is arbitrary and so by computing this minimum we get an expression for $u^+(\zt)$, for every $\zt\in[0,T]$.

It turns out that  $\frac{\ZD}{\ZD t}W(\zt  ;X _{\zt  })+u^*(\zt)u(\zt)$ is the sum of several terms. Some of them do not depend on $u$ and the minimization concerns solely the terms which depends on $u$. We get (we recall that $P_0(\zt)$ is selfadjoint)
\begin{align}
\nonumber u^+(\zt)&={\rm arg\, min} \left \{
u^*B^*P_0(\zt)\hat x +u^*B^* \int_0^\zt P_1(\zt,s) \tilde x(s)\ZD s+\right.\\
\ZLA{FunzioDAMINIMperilCONTROTTIMO}&\left.+\hat x^*P_0(\zt)Bu+\left ( \int_0^\zt \tilde x^*(s)P_1^*(\zt,s)\ZD s\right )Bu+u^*u
\right \}\,.
\end{align}
The minimization gives  
\begin{equation}
\ZLA{eq:FeedbackFORMuOTTIMOpre}
u^+(\zt)=-B^*\left [P_0(\zt) \hat x +\int_0^\zt P_1(\zt,s)\tilde x (s)\ZD s\right ]\,.
\end{equation}
If the system is solved up to time $t$ along an optimal evolution (so that $x^+(t)$ is equal to $\tilde x(t)$ when $t<\zt$ and it is the solution which corresponds to the  optimal control for larger times) we have
\[
\ZLA{eq:FeedbackFORMuOTTIMO}
u^+(t)=-B^*\left [P_0(\zt)   x^+(t) +\int_0^\zt P_1(t,s)  x^+ (s)\ZD s\right ] 
\]
and this is the feedback form of the optimal control (compare~\cite{Pritch}).

We repalce~(\ref{eq:FeedbackFORMuOTTIMOpre}) in the brace in~(\ref{FunzioDAMINIMperilCONTROTTIMO}) and we see that the minimum is
\begin{align}
\nonumber&-\hat x^*P_0(\zt)BB^*P_0(\zt)\hat x-\hat x^*P_0(\zt)BB^*\int_0^\zt P_1(\zt,\xi)\tilde x(\xi)\ZD\xi-\\
\ZLA{eq:IlMINIMOallOOTTIMMMO}&-\left (\int_0^\zt \tilde x^*(r)P_1^*(\zt,r)\ZD r\right )BB^*P_0(\zt)\hat x-\int_0^\zt\int_0^\zt \tilde x(r)P_1(\zt,r)BB^*P_1(\zt,\xi)\tilde x(\xi)\ZD\xi\ZD r\,.
\end{align}

Now we compute the derivative of the function $ \zt\mapsto W(\zt  ;X _{\zt  })$ along an optimal evolution and we consider its limit for $t\to \zt+$. We insert this quantity in~(\ref{eq:PrimaFormINEqEQ}), which is an equality since we are computing the limit along an optimal evolution. We take into account that the terms which contains $u$ sum up to the expression~(\ref{eq:IlMINIMOallOOTTIMMMO}) and we get the following equality. In this equality, a superimposed dot denotes derivative with respect to the variable $\zt$:
\[
\dot P_0(\zt)=\frac{\ZD}{\ZD\zt} P_0(\zt)\,,\quad 
\dot P_1(\zt,\xi)=\frac{\partial}{\partial \zt} P_1(\zt,\xi)\,,\qquad \dot K(\zt,\xi,r)=
\frac{\partial}{\partial \zt} K(\zt,\xi,r)\,.
\]
The equality is:
\begin{align*}
 &-\hat x^*P_0(\zt)BB^*P_0(\zt)\hat x-\hat x^*P_0(\zt)BB^*\int_0^\zt P_1(\zt,\xi)\tilde x(\xi)\ZD\xi-\\
&-\left (\int_0^\zt \tilde x^*(r)P_1^*(\zt,r)\ZD r\right )BB^*P_0(\zt)\hat x-\int_0^\zt\int_0^\zt \tilde x(r)P_1(\zt,r)BB^*P_1(\zt,\xi)\tilde x(\xi)\ZD\xi\ZD r+\\
&+\left (\int_0^\zt \tilde x^*(r)N^*(\zt-r)\ZD s\right )P_0(\zt)\hat x+\hat x^*\dot P_0(\zt)\hat x+\hat x^*\int_0^\zt N(\zt-\xi)\tilde x(\xi)\ZD\xi+ 
 \hat x^*P_1(\zt,\zt)\hat x+\\
 &+\hat x^*P_1^*(\zt,\zt)\hat x+\left ( \int_0^\zt \tilde x^*(r)N^*( \zt-r)\ZD r\right )\int_0^\zt P_1(\zt,s)\tilde x(s)\ZD s+\hat x^*\int_0^\zt \dot P_1(\zt,\xi)\tilde x(\xi)\ZD \xi+\\
 &+\left (\int_0^\zt \tilde x^*(r)\dot P_1^*(\zt,r)\ZD r\right )\hat x +\left ( \int_0^\zt \tilde x^*(r)P_1(\zt,r)\ZD r\right )\left ( \int_0^\zt N(\zt-\xi)\tilde x(\xi)\ZD \xi\right )+\\
& +\left (\int_0^\zt \tilde x^*(r)K(\zt,\zt,r)\ZD r\right )\hat x+\hat x^*\int_0^\zt K(\zt,\xi,\zt)\tilde x(\xi)\ZD\xi+\\
&+\int_0^\zt \tilde x^*(r)\int_0^\zt \dot K(\zt,\xi,r)\tilde x(\xi)\ZD\xi\ZD r+\hat x^*Q\hat x=0
\end{align*}

The vector $\hat x$ and the function $\tilde x(\cdot)$ are arbitrary. So, we first impose $\tilde x(\cdot)=0$ and $\hat x$ arbitrary, then the converse and finally both nonzero arbitrary. We find that the three matrix functions $P_0(\zt)$, $P_1(\zt,r)$, $K(\zt,\xi,r)$ solve the following system of differential equations in the arbitrary variable $\zt$. The variables $r$ and $\xi$ belong to $[0,\zt]$ for every $\zt\in [0,T]$.

\begin{equation}\ZLA{Eq:RICCATI}
\begin{array}{ll}
&\displaystyle 	P_0'(\zt)-P_0(\zt)B^*BP_0(\zt)+Q(\zt)+P_1(\zt,\zt)+P_1^*(\zt,\zt)=0\\[2mm]
&\displaystyle \frac{\partial}{\partial \zt}P_{1}(\zt,\xi)  -P_0(\zt)BB^*P_1(\zt,\xi) +P_0(\zt)N(\zt-\xi) 
   +K(\zt,\xi,\zt)=0\\[2mm]
&\displaystyle \frac{\partial}{\partial \zt}K(\zt,\xi,r)-P_1^*(\zt,r)BB^*P_1(\zt,\xi)+\\
 &\displaystyle{~}\quad +P^*_1(\zt,r)N(\zt-\xi)
 +N^*(\zt-r)P_1(\zt,\xi) =0\\[2mm]
&\displaystyle P_0(T)=Q_0\,,\qquad P_1(T,\xi)=0\,,\qquad K(T,\xi,r)=0
\end{array}
\end{equation}

The final conditions are obtained by noting that when $\zt=T$ i.e. with $X_T=(\hat x,\tilde x_T(\cdot))$ arbitrary in $M^2_T=\zzr^n\times L^2(0,T;\zzr^n)$, the expression $W(T,X_T)$ in~(\ref{eq:ExprreDELLAformAquaDRATI}) is equal to $J_T(X_T;u)=\hat x^*Q_0\hat x$ for every $X_T$.

\emph{This is the Riccati differential equation of our optimization problem.}

\begin{Remark}
{\rm 
We note the following facts:
\begin{itemize}
\item
We take into account the fact that $P_0$ is selfadjoint and $K^*(\zt,\xi,\zt)=K(\zt,\zt,\xi)$. We compute the adjoint of the second line in~(\ref{Eq:RICCATI}) and we find:
\[
 \frac{\partial}{\partial \zt}P^*_{1}(\zt,r)  -P_1^*(\zt,r)BB^*P_0(\zt)
 + N^*(\zt-r)P_0(\zt)  +K(\zt,\zt,r) =0 \,. 
\]
\item
The form of the Riccati differential equations we derived for the Volterra integral equation~(\ref{eq:Volte}) has to be compared with the Riccati differential equation `` in decoupled form'' which was once fashionable in the study of the quadratic regulator problem for systems with finite delays, see~\cite{Ross}.\zdia
\end{itemize}
}
\end{Remark}

\end{document}